# НОВЫЙ АЛГОРИТМ ПРОГОНКИ РЕШЕНИЯ НЕПРЕРЫВНОЙ ЛИНЕЙНО КВАДРАТИЧНОЙ ЗАДАЧИ ОПТИМИЗАЦИИ С НЕРАЗДЕЛЕННЫМИ КРАЕВЫМИ УСЛОВИЯМИ

Ф.А.Алиев, М.М.Муталлимов, И.А.Магеррамов,

Н.А.Исмайлов, С.М.Мирсаабов

Институт Прикладной Математики, БГУ, Баку, Азербайджан

e-mail: f_aliev@yahoo.com

**Резюме.** Приводится новый алгоритм прогонки решения линейно- квадратичной задачи (ЛКЗ) оптимизации с неразделенными краевыми условиями в непрерывном случае. Используя свойства $J$ - симметричности соответствующей Гамильтоновой матрицы, уравнений Эйлера-Лагранжа, показывается, что линейные алгебраические уравнения для определения недостающих начальных данных решаемой системы имеют симметричную главную матрицу. Результаты иллюстрируются на примере ЛКЗ оптимизации (стационарный случай) с минимальными управляющими воздействиями.

**Ключевые слова:** оптимизация, метод прогонки, дифференциальные уравнение Риккати, линейные системы алгебраические уравнение.

**AMS Subject Classification:** 49J15, 49J35.

### Введение

Как известно [9,14,21] для решения линейно квадратичной задачи оптимизации с неразделенными краевыми условиями в непрерывном случае [7,8,12,21] существуют разные методы – метод, повышающий размерности исходной системы [19,20], метод прогонки [10,18], метод Мощинского [4,6] и т.д. Однако, каждый из этих методов сталкивается с трудностями [1,5] при определенных случаях: например, обобщение метода [10] к многоточечному случаю с неразделенными краевыми условиями, переход через узловые точки сталкивается с серьезными трудностями из-за не единственности переходов и т.п. Поэтому, в данном случае предлагается новый метод, не требующий решения матричных уравнений Риккати, линейных матричных уравнений и др. Далее, для построения соответствующих фундаментальных матриц приводится алгоритм, требующий решения матричных дифференциальных уравнений гораздо меньшей размерности. Результаты иллюстрируются на примере стационарной ЛКЗ оптимизации с минимальным управляющим воздействием.

## 2. Постановка задачи

Пусть движение описывается системой линейных дифференциальных уравнений

$$\dot{x} = F(t)x(t) + G(t)u(t) \qquad (1)$$

с неразделенными краевыми условиями

$$\Phi_1 x(t_0) - \Phi_2 x(\tau) = q \qquad (2)$$

где $x$ - $n$ мерный фазовый вектор, $u(t)$ - $m$ мерный вектор управляющих воздействий, $F(t)$, $G(t)$ - известные кусочно-непрерывные функции-матрицы $n \times n$, $n \times m$ размерности, $\Phi_1$, $\Phi_2$ -постоянные матрицы $k \times n$ размерности, $q$ - постоянный $k \times 1$ мерный вектор, время $\tau$ заданные числа. Предполагается, что пара $(F(t), G(t))$, управляемая в каждой точке отрезка времени $(t_0, \tau)$ [13], а $[\Phi_1, \Phi_2, q]$ матрица, которая удовлетворяет условию Кронеккера-Капелли [2,3,11].

Требуется найти управляющие воздействия $u(t)$ так, чтобы с соответствующим $x(t)$ из (1), (2) минимизировали квадратичный критерий качества

$$J = \frac{1}{2}\int_0^\tau \left(x'(t)R(t)x(t) + u'(t)C(t)u(t)\right)dt \qquad (3)$$

Здесь $R(t) = R'(t) \geq 0$, $C(t) = C'(t) > 0$, соответственно $n \times n$, $m \times m$ размерные, тоже кусочно-непрерывные функции-матрицы, а штрих означает операцию транспонирования.

## 3. Метод, повышающий размерности исходной системы

Как известно [15], решение задачи (1)-(3) сводится к нахождению решения следующей системы уравнений Эйлера-Лагранжа

$$\begin{vmatrix} \dot{x} \\ \dot{\lambda} \end{vmatrix} = \begin{bmatrix} F(t) & -G(t)C^{-1}(t)G'(t) \\ Q(t) & -F'(t) \end{bmatrix} \begin{bmatrix} x \\ \lambda \end{bmatrix} = H(t)\begin{bmatrix} x(t) \\ \lambda(t) \end{bmatrix} \qquad (4)$$

с дополнительными краевыми условиями

$$\begin{aligned}\lambda(t_0) &= -\Phi_1'\nu \\ \lambda(\tau) &= -\Phi_2'\nu\end{aligned} \qquad (5)$$

где $\lambda(t)$ и $\nu$ соответствующие сопряженные векторы Лагранжа оптимизационной задачи (1) - (3).

Пусть $\Phi(t, t_0)$ является фундаментальной матрицей системы (4), т.е.

$$\dot{\Phi}(t, t_0) = H(t)\Phi(t, t_0), \quad \Phi(t_0, t_0) = E \qquad (6)$$

где $E$ - $n \times n$ единичная матрица. Обозначая

$$\Phi(t, t_0) = \begin{bmatrix} \Phi_{11}(t, t_0) & \Phi_{12}(t, t_0) \\ \Phi_{21}(t, t_0) & \Phi_{22}(t, t_0) \end{bmatrix} \qquad (7)$$

решение системы (4) представим в виде

$$x(t) = \Phi_{11}(t,t_0)\,x(t_0) + \Phi_{12}(t,t_0)\,\lambda(t_0) \\ \lambda(t) = \Phi_{21}(t,t_0)\,x(t_0) + \Phi_{22}(t,t_0)\,\lambda(t_0) \quad (8)$$

Пусть из (8) $\Phi_{22}^{-1}(t,t_0)$ существует, тогда из второго уравнения находим $\lambda(t_0)$ и, подставив в первое уравнение, получим

$$x(t) = \left(\Phi_{11}(t,t_0) - \Phi_{12}(t,t_0)\Phi_{22}^{-1}(t,t_0)\Phi_{21}(t,t_0)\right)x(t_0) + \Phi_{12}(t,t_0)\Phi_{22}^{-1}(t,t_0)\lambda(t) \quad (9)$$

$$\lambda(t_0) = -\Phi_{22}^{-1}(t,t_0)\Phi_{21}(t,t_0)x(t_0) + \Phi_{22}^{-1}(t,t_0)\lambda(t).$$

При $t = \tau$, решив системы линейно алгебраических уравнений (9), (5), (2), находим $x(0)$, $\lambda(0)$, $x(\tau)$, $\lambda(\tau)$ и $\nu$. Далее из (8) находим текущие значения $x(t)$, $\lambda(t)$ из (8). А управляющие воздействие определяются из соотношений.

$$u(t) = -C^{-1}(t)G(t)\lambda(t). \quad (10)$$

Такой подход является методом, повышающий размерность исходной системы (1). Поэтому, попытаемся уменьшить размерность в процессе решения задачи (1)-(3).

### 4. Новый метод прогонки

Для уменьшения вычислений используем соотношение (5) при $t = \tau$ в уравнениях (9), т.е. напишем (9) при $t = \tau$

$$x(\tau) = \left(\Phi_{11}(\tau,t_0) - \Phi_{12}(\tau,t_0)\Phi_{22}^{-1}(\tau,t_0)\Phi_{12}(\tau,t_0)\right)x(t_0) + \\ + \Phi_{12}(\tau,t_0)\Phi_{22}^{-1}(\tau,t_0)\lambda(\tau) \\ \lambda(t_0) = -\Phi_{22}^{-1}(\tau,t_0)\Phi_{21}(\tau,t_0)x(t_0) + \Phi_{22}^{-1}(\tau,t_0)\lambda(\tau) \quad (11)$$

Далее, учитывая (5) в (10), имеем

$$x(\tau) = \left(\Phi_{11}(\tau,t_0) - \Phi_{12}(\tau,t_0)\Phi_{22}^{-1}(\tau,t_0)\Phi_{21}(\tau,t_0)\right)x(t_0) - \Phi_{12}(\tau,t_0)\Phi_{22}^{-1}(\tau,t_0)\Phi_2'\nu \\ -\Phi_1'\nu = -\Phi_{22}^{-1}(\tau,t_0)\Phi_{21}(\tau,t_0)x(t_0) - \Phi_{22}^{-1}(\tau,t_0)\Phi_2'\nu \quad (12)$$

Подставляя из (12) $x(\tau)$ в (2) и объединяя со вторым уравнением (12) для определения $x(t_0)$ и $\nu$ имеем следующую систему алгебраических уравнений

$$\Phi_{22}^{-1}(\tau,t_0)\Phi_{21}(\tau,t_0)x(t_0) + \left(-\Phi_{22}^{-1}(\tau,t_0)\Phi_2' + \Phi_1'\right)\nu = 0 \\ \left[-\Phi_1 - \Phi_2(\Phi_{11}(\tau,t_0) - \Phi_{12}(\tau,t_0)\Phi_{22}^{-1}(\tau,t_0)\Phi_{21}(\tau,t_0))\right]x(t_0) + \\ + \Phi_2\Phi_{12}(\tau,t_0)\Phi_{22}^{-1}(\tau,t_0)\Phi_2'\nu = q \quad (13)$$

Теперь докажем симметричность главной матрицы (13), т.е. из (13)

$$\Phi_{22}^{-1}(\tau,t_0)\Phi_{21}(\tau,t_0) = \Phi_{21}'(\tau,t_0)\left(\Phi_{22}^{-1}(\tau,t_0)\right)' \quad (14)$$

$$\Phi_{22}^{-1}(\tau,t_0)\Phi_2' - \Phi_1' = \left[\Phi_1 - \Phi_2\Phi_{11}(\tau,t_0) - \Phi_{12}(\tau,t_0)\Phi_{22}^{-1}(\tau,t_0)\Phi_{21}'(\tau,t_0)\right] \quad (15)$$

$$\Phi_2\Phi_{12}(\tau,t_0)\Phi_{22}^{-1}(\tau,t_0)\Phi_2' = \Phi_2\left(\Phi_{22}^{-1}(\tau,t_0)\right)'\Phi_{12}'(\tau,t_0)\Phi_2' \quad (16)$$

Используя симметричность системы (4), из [14]

$$JHJ' = -H, \qquad J = \begin{bmatrix} 0 & E \\ -E & 0 \end{bmatrix} \qquad (17)$$

легко доказывается, что при помощи соотношений из $J\Phi(t)J'\Phi'(t) = E$ получаем

$$E = \begin{bmatrix} \Phi_{22}(\tau,t_0)\Phi'_{11}(\tau,t_0) - \Phi_{21}(\tau,t_0)\Phi'_{12}(\tau,t_0) & \Phi_{22}(\tau,t_0)\Phi'_{21}(\tau,t_0) - \Phi_{21}(\tau,t_0)\Phi'_{22}(\tau,t_0) \\ \Phi_{11}(\tau,t_0)\Phi'_{12}(\tau,t_0) - \Phi_{12}(\tau,t_0)\Phi'_{11}(\tau,t_0) & \Phi_{11}(\tau,t_0)\Phi'_{22}(\tau,t_0) - \Phi_{12}(\tau,t_0)\Phi'_{21}(\tau,t_0) \end{bmatrix} \qquad (18)$$

Из условий (18) верхнюю вне диагональную подматрицу напишем как

$$\Phi_{22}(\tau,t_0)\Phi'_{21}(\tau,t_0) - \Phi_{21}(\tau,t_0)\Phi'_{22}(\tau,t_0) = 0. \qquad (19)$$

Умножая с правой стороны (19) на $\left(\Phi'_{22}(\tau,t_0)\right)^{-1}$, а с левой стороны на $\left(\Phi^{-1}_{22}(\tau,t_0)\right)$, получим соотношение (14). Таким образом доказана следующая Лемма.

**Лемма 1.** Из системы линейных матричных алгебраических уравнений (13) подматрицы $\Phi^{-1}_{22}(\tau,t_0)\,\Phi_{21}(\tau,t_0)$ - симметричные подматрицы.

Теперь покажем, что из (13) матрица $\Phi_1(\tau,t_0) \cdot \Phi^{-1}_{22}(\tau,t_0)$ -симметричная, т.е.

$$\Phi_1(\tau,t_0) \cdot \Phi^{-1}_{22}(\tau,t_0) = \Phi^{-1\,\prime}_{22}(\tau,t_0)\Phi'_{12}(\tau,t_0) \qquad (20)$$

Тогда нижняя диагональная подматрица в (13) тоже будет симметричной. На самом деле

$$\left(\Phi_2\Phi_{12}(\tau,t_0)\Phi^{-1}_{22}(\tau,t_0)\Phi'_2\right)' = \Phi_2\left(\Phi_{12}(\tau,t_0)\Phi^{-1}_{22}(\tau,t_0)\right)'\Phi'_2 =$$
$$= \Phi_2\left(\Phi^{-1}_{22}(\tau,t_0)\right)'\Phi'_{12}(\tau,t_0)\Phi'_2 = \Phi_2\Phi_{12}(\tau,t_0)\Phi^{-1}_{22}(\tau,t_0)\Phi'_2$$

Теперь докажем равенство (20). Из (18) следует, что

$$\Phi_{22}(\tau,t_0)\Phi'_{11}(\tau,t_0) - \Phi_{21}(\tau,t_0)\Phi'_{12}(\tau,t_0) = E \qquad (21)$$

$$\Phi_{22}(\tau,t_0)\Phi'_{21}(\tau,t_0) - \Phi_{21}(\tau,t_0)\Phi'_{22}(\tau,t_0) = 0 \qquad (22)$$

$$\Phi_{11}(\tau,t_0)\Phi'_{12}(\tau,t_0) - \Phi_{12}(\tau,t_0)\Phi'_{11}(\tau,t_0) = 0 \qquad (23)$$

$$\Phi_{11}(\tau,t_0)\Phi'_{22}(\tau,t_0) - \Phi_{12}(\tau,t_0)\Phi'_{21}(\tau,t_0) = E \qquad (24)$$

Умножая с левой стороны (22) на $\Phi^{-1}_{22}(\tau,t_0)$ получим

$$\Phi'_{21}(\tau,t_0) = \Phi^{-1}_{22}(\tau,t_0)\Phi_{21}(\tau,t_0)\Phi'_{22}(\tau,t_0)$$

Далее, транспонируя обе части полученного последнего выражения, будем иметь

$$\Phi_{21}(\tau,t_0) = \Phi_{22}(\tau,t_0)\Phi'_{21}(\tau,t_0)\left(\Phi^{-1}_{22}(\tau,t_0)\right)'$$

Умножим полученное равенство с левой стороны на $\Phi^{-1}_{22}(\tau,t_0)$

$$\Phi^{-1}_{22}(\tau,t_0)\Phi_{21}(\tau,t_0) = \Phi'_{21}(\tau,t_0)\left(\Phi^{-1}_{22}(\tau,t_0)\right)'. \qquad (25)$$

Сейчас, умножая (21) с левой стороны на $\Phi^{-1}_{22}(\tau,t_0)$, получим

$$\Phi'_{11}(\tau,t_0) = \Phi^{-1}_{22}(\tau,t_0) + \Phi^{-1}_{22}(\tau,t_0)\Phi_{21}(\tau,t_0)\Phi'_{12}(\tau,t_0).$$

Подставляя полученные выражения $\Phi'_{11}(\tau,t_0)$ в (23) имеем

$$\Phi_{11}(\tau,t_0)\Phi'_{12}(\tau,t_0) = \Phi_{12}(\tau,t_0)\Phi'_{11}(\tau,t_0) =$$
$$= \Phi_{12}(\tau,t_0)\left[\Phi_{22}^{-1}(\tau,t_0) + \Phi_{22}^{-1}(\tau,t_0)\Phi_{21}(\tau,t_0)\Phi'_{12}(\tau,t_0)\right] =$$
$$= \Phi_{12}(\tau,t_0)\Phi_{22}^{-1}(\tau,t_0) + \Phi_{12}(\tau,t_0)\Phi_{22}^{-1}(\tau,t_0)\Phi_{21}(\tau,t_0)\Phi'_{12}(\tau,t_0).$$

Отсюда
$$\Phi_{12}(\tau,t_0)\Phi_{22}^{-1}(\tau,t_0) = \left[\Phi_{11}(\tau,t_0) - \Phi_{12}(\tau,t_0)\Phi_{22}^{-1}(\tau,t_0)\Phi_{21}(\tau,t_0)\right]\Phi'_{12}(\tau,t_0).$$

Учитывая здесь (25), получим
$$\Phi_{12}(\tau,t_0)\Phi_{22}^{-1}(\tau,t_0) = \left[\Phi_{11}(\tau,t_0) - \Phi_{12}(\tau,t_0)\Phi'_{21}(\tau,t_0)\left(\Phi_{22}^{-1}(\tau,t_0)\right)'\right]\Phi'_{12}(\tau,t_0), \qquad (26)$$

и определим из (24) $\Phi_{11}(\tau,t_0)$
$$\Phi_{11}(\tau,t_0) = \left(\Phi'_{22}(\tau,t_0)\right)^{-1} + \Phi_{12}(\tau,t_0)\Phi'_{21}(\tau,t_0)\left(\Phi'_{22}(\tau,t_0)\right)^{-1}$$

значение которого подставив в (26), имеем
$$\Phi_{12}(\tau,t_0)\Phi_{22}^{-1}(\tau,t_0) =$$
$$= \left[\left(\Phi'_{22}(\tau,t_0)\right)^{-1} + \Phi_{12}(\tau,t_0)\Phi'_{21}(\tau,t_0)\left(\Phi'_{22}(\tau,t_0)\right)^{-1} - \Phi_{12}(\tau,t_0)\Phi'_{21}(\tau,t_0)\left(\Phi_{22}^{-1}(\tau,t_0)\right)^1\right]\Phi'_{12}(\tau,t_0) = .$$
$$= \left(\Phi_{22}^{-1}(\tau,t_0)\right)'\Phi'_{12}(\tau,t_0)$$

Таким образом, доказана следующая

**Лемма 2.** В системе линейных матричных алгебраических уравнений (13) матрица $\Phi_2$ $\Phi_{12}(\tau,t_0)$ $\Phi_{22}^{-1}(\tau,t_0)$ $\Phi'_2$ является симметричной.

Теперь докажем, что вне диагональные матрицы уравнения (13) являются симметричными, т.е.

$$\left[\Phi_1 - \Phi_2\left(\Phi_{11}(\tau,t_0) - \Phi_{12}(\tau,t_0)\Phi_{22}^{-1}(\tau,t_0)\Phi_{21}(\tau,t_0)\right)\right] = \left(-\Phi_{22}^{-1}(\tau,t_0)\Phi'_2 + \Phi'_1\right)'. \qquad (27)$$

Для выполнения (26) достаточно показать, что
$$\Phi'^{-1}_{22}(\tau,t_0) = \Phi_{11}(\tau,t_0) - \Phi_{12}(\tau,t_0)\Phi_{22}^{-1}(\tau,t_0)\Phi_{21}(\tau,t_0). \qquad (28)$$

Действительно из (28) и (21)-(24) можно показать, что
$$\Phi_{11}(\tau,t_0) - \Phi_{12}(\tau,t_0)\Phi_{22}^{-1}(\tau,t_0)\Phi_{21}(\tau,t_0) = \Phi_{11}(\tau,t_0) - \Phi'^{-1}_{22}(\tau,t_0)\Phi'_{12}(\tau,t_0)\Phi_{21}(\tau,t_0) =$$
$$= \left(\Phi'_{22}(\tau,t_0)\right)^{-1} + \Phi_{12}(\tau,t_0)\Phi'_{21}(\tau,t_0)\Phi'^{-1}_{22}(\tau,t_0) - \Phi'^{-1}_{22}(\tau,t_0)\Phi'_{12}(\tau,t_0)\Phi_{21}(\tau,t_0) =$$
$$= \Phi'^{-1}_{22}(\tau,t_0) + \Phi_{12}(\tau,t_0)\Phi'_{21}(\tau,t_0)\Phi_{22}^{-1}(\tau,t_0) - \Phi_{12}(\tau,t_0)\Phi_{22}^{-1}(\tau,t_0)\Phi_{21}(\tau,t_0) =$$
$$= \Phi'^{-1}_{22}(\tau,t_0) + \Phi_{12}(\tau,t_0)\Phi'_{21}(\tau,t_0)\Phi'^{-1}_{22}(\tau,t_0) - \Phi_{12}(\tau,t_0)\Phi'_{21}(\tau,t_0)\Phi'^{-1}_{22}(\tau,t_0) = \Phi'^{-1}_{22}(\tau,t_0),$$

т.е. удовлетворяется условие симметричности (28).

Таким образом доказана следующая

**Лемма 3.** В системе линейных матричных уравнений (13) удовлетворяется условие симметричности (27).

Результаты Лемм 1-3 позволяют нам делать вывод, что главная матрица системы алгебраических уравнений (13) является симметричной. Для этого сначала уравнение запишем в матричном виде, т.е., обозначая

$$D = \begin{bmatrix} -\Phi_{22}^{-1}(\tau,t_0)\,\Phi_{21}(\tau,t_0) & \left(-\Phi_{22}^{-1}(\tau,t_0)\Phi_2' + \Phi_1'\right) \\ \Phi_1 - \Phi_2\,(\Phi_{11}(\tau,t_0) - \Phi_{12}(\tau,t_0)\Phi_{22}^{-1}(\tau,t_0)\Phi_{21}(\tau,t_0)) & \Phi_2\,\Phi_{12}(\tau,t_0)\Phi_{22}^{-1}(\tau,t_0)\Phi_2' \end{bmatrix}, \quad (29)$$

$$\widetilde{W} = \begin{bmatrix} 0 \\ q \end{bmatrix}$$

напишем (13) в следующем более компактном виде

$$D\begin{bmatrix} x(t_0) \\ \nu \end{bmatrix} = \widetilde{W} \quad . \tag{30}$$

Имеем следующую теорему.

**Теорема.** Матрица $D$ из алгебраических уравнений (13), (30) является симметричной.

Поскольку в выражениях (9) и (10) для определения $x(t)$ и $u(t)$ содержится $\lambda(t)$, попытаемся получить выражения для $x(t)$ и $u(t)$ без использования $\lambda(t)$. Для этого, выражение $\lambda(t)$ из (8) поставим в первое выражение (9), а для определения $x(t)$ и $u(t)$ учитываем $\lambda(t_0) = -\Phi_1'\nu$. Тогда для $x(t)$ получим

$$x(t) = \Phi_{11}(\tau,t_0)x(t_0) - \Phi_{12}(\tau,t_0)\Phi_1'\,\nu \tag{31}$$

$$u(t) = -C^{-1}(t)G(t)\Phi_{21}(\tau,t_0)x(t_0) + C^{-1}(t)G(t)\Phi_{22}(\tau,t_0)\Phi_1'\nu \tag{32}$$

Таким образом имеется следующий

Алгоритм.

1. Формируются матрицы $F(t)$, $G(t)$, $\Phi_1$, $\Phi_2$, $R(t)$, $C(t)$ из (1)-(3)

2. Формируется матрица $H(t)$ из (4)

3. Решается матричное уравнение (6) и находится $\Phi(t,t_0)$ из (7)

4. Восстанавливается матрица $D = D'$ и вектор $\widetilde{W}$ согласно (29), и решается матричное алгебраическое уравнение (30), находятся $x(t_0)$, $\nu$.

5. По формуле (31), (32) вычисляются $x(t)$ и $u(t)$ через $x(t_0)$ и $\nu$.

Теперь приведем формулы для $u(t)$ через фазовые координаты $x(t)$. Для этого напишем в (8) первое соотношение из (5)

$$x(t) = \Phi_{11}(t,t_o)x(t_0) + \Phi_{12}(t,t_0)\Phi_1\nu$$

$$\lambda(t) = \Phi_{21}(t,t_o)x(t_0) - \Phi_{22}(t,t_0)\Phi_1\nu \tag{33}$$

Тогда при предположении существования $\Phi_{11}^{-1}(t,t_0)$ находим $x(t_0)$ из первого соотношения (33) в виде

$$x(t_0) = \Phi_{11}^{-1}(t,t_o)x(t_0) + \Phi_{22}(t,t_0)\Phi_1'\nu \qquad (34)$$

Подставив (34) в последнее соотношение (33) получим

$$\lambda(t) = \Phi_{21}(t,t_0)\Phi_{11}^{-1}(t,t_0)x(t) + (\Phi_{21}(t,t_o)\Phi_{11}^{-1}(t,t_0)\Phi_{12}(t,t_o) - \Phi_{22}(t,t_0)\Phi_1'\nu \qquad (35)$$

который является выражением $\lambda(t)$ через фазовые состояния $x(t)$. Подставив выражение $\lambda(t)$ из (35) в (10) получим для управляющих воздействий следующее выражение

$$\begin{aligned}u(t) = -C^{-1}(t)G(t)[&\Phi_{21}(t,t_0)\Phi_{11}^{-1}(t,t_o)x(t) + \\ + (&\Phi_{21}(t,t_0)\Phi_{11}^{-1}(t,t_o)\Phi_{12}(t,t_0) - \Phi_{22}(t,t_0))\Phi_1'\nu]\end{aligned} \qquad (36)$$

которое является обратной связью для решения задачи (1)-(3).
Подставляя (36) в (1) для нахождения $x(t)$ имеем следующую систему дифференциальных уравнений

$$\begin{aligned}\dot{x}(t) = \Big(F(t) - G(t)C^{-1}(t)G(t)\Phi_{21}(t,t_0)\Phi_{11}^{-1}(t,t_o)\Big)x(t) + \\ + G(t)C^{-1}(t)G(t)\Big(\Phi_{21}(t,t_0)\Phi_{11}^{-1}(t,t_o)\Phi_{12}(t,t_0) - \Phi_{22}(t,t_0)\Big)\Phi_1'\nu\end{aligned} \qquad (37)$$

с начальными условиями $x(t_0)$, найденными из системы линейных алгебраических уравнений (30). В этом случае в вышеизложенном алгоритме пункт 5 заменяется следующим пунктом

5. По формуле (36) определяется $u(t)$, а $x(t)$ находится из решения дифференциальных уравнений (37) с начальным условием $x(t_0)$.

**5. Построение фундаментальных матриц (6)**

Таким образом, решение задачи (1)-(3) сводится к нахождению фундаментальной матрицы системы (4), т.е. к решению задачи (6). Однако, решение задачи (6) в общем случае представляется возможным с трудностями (из-за размерности плохой обусловленности матриц $H$ и т.д.). Поэтому для нахождения $\Phi(\tau,t_0)$ воспользуемся методом Захар-Иткина [16]. Действительно из [16,17]

$$\begin{cases} x(t) = \psi x(t_0) - W\lambda(t) \\ \lambda(t_0) = Vx(t_0) + \psi'\lambda(t), \end{cases} \quad (38)$$

где $\psi$, $W$, $V$ матрицы, которые удовлетворяют следующим дифференциальным уравнениям и начальным условиям:

$$\begin{cases} \psi' = (F + WR)\psi, & \psi(t_0) = E \\ \dot{W} = FW + WF' + WRW - GQ^{-1}G', & W(t_0) = 0 \\ \dot{V} = \psi'R\psi, & V(t_0) = 0. \end{cases} \quad (39)$$

Приводя некоторые преобразования из (33) получим

$$\begin{bmatrix} x(t) \\ \lambda(t) \end{bmatrix} = \begin{bmatrix} E & W \\ 0 & -\psi' \end{bmatrix}^{-1} \cdot \begin{bmatrix} \psi & 0 \\ V & -E \end{bmatrix} \begin{bmatrix} x(t_0) \\ \lambda(t_0) \end{bmatrix} = \begin{bmatrix} E & W(\psi')^{-1} \\ 0 & -(\psi')^{-1} \end{bmatrix} \cdot \begin{bmatrix} \psi & 0 \\ V & -E \end{bmatrix} \begin{bmatrix} x(t_0) \\ \lambda(t_0) \end{bmatrix} =$$
$$= \begin{bmatrix} \psi + W(\psi')^{-1}V & -W(\psi')^{-1} \\ -(\psi')^{-1}V & (\psi')^{-1} \end{bmatrix} \begin{bmatrix} x(t_0) \\ \lambda(t_0) \end{bmatrix}. \quad (40)$$

Отсюда следует, что фундаментальная матрица (7), образованная через $\Phi_{11}(t,t_0)$, $\Phi_{12}(t,t_0)$, $\Phi_{21}(t,t_0)$, $\Phi_{22}(t,t_0)$ определяется как

$$\Phi_{11}(t,t_0) = \psi(t,t_0) + W(t,t_0) \cdot (\psi'(t,t_0))^{-1} V(t,t_0),$$
$$\Phi_{12}(t,t_0) = -W(t,t_0) \cdot (\psi'(t,t_0))^{-1},$$
$$\Phi_{21}(t,t_0) = -(\psi'(t,t_0))^{-1} V(t,t_0), \quad (41)$$
$$\Phi_{22}(t,t_0) = (\psi'(t,t_0))^{-1}.$$

Приведем следующий алгоритм программ на основе вышесказанных результатов. Иллюстрируем вышеизложенный алгоритм на следующем примере.

**Пример**

Пусть в (1)-(3) $F(t) = F$, $G(t) = G$, $R(t) = 0, C(t) = C$ постоянные матрицы с соответствующими размерностями. Для начала восстанавливаем фундаментальную матрицу $\Phi(t,t_0)$ через матричные дифференциальные уравнения (39) и формулы (40)-(41). На самом деле, в данном случае (39) упрощается и переходит к более простому виду

$$\dot{\psi} = F\psi, \quad \psi(t_0) = E$$
$$\dot{W} = FW + WF' - GQ^{-1}G', \quad W(t_0) = 0 \quad (42)$$
$$\dot{V} = 0 \quad V(t_0) = 0$$

где решение (42) в аналитическом виде представляется как

$$\psi(t,t_0) = e^{F(t-t_0)},$$

$$W(t,t_0) = W_1 - e^{F(t-t_0)} W_1 e^{F'(t-t_0)} \quad (43)$$

$$V(t,t_0) = 0$$

где $W_1$ является решением следующего матричного уравнения Ляпунова

$$-FW_1 - W_1 F' + GQ^{-1}G' = 0 \quad (44)$$

Отметим, что уравнение Ляпунова имеет решение в следующем виде [13]

$$W_1 = -\int_0^\infty e^{Ft} G_1 Q^{-1} C_1' e^{F't} dt \quad (45)$$

Теперь формируем фундаментальную матрицу $\Phi(t,t_0)$ через решения уравнений (42)

$$\Phi_{11}(t,t_0) = e^{F(t-t_0)}$$

$$\Phi_{12}(t,t_0) = -[W_1 - e^{F(t-t_0)} W_1 e^{F'(t-t_0)}] e^{-F(t-t_0)} \quad (46)$$

$$\Phi_{21}(t,t_0) = 0$$

$$\Phi_{22}(t,t_0) = e^{-F(t-t_0)}$$

Подставляя (46) в (36) имеем для управления следующее выражение

$$u(t) = C^{-1}(t) G(t) e^{-F(t-t_0)} \Phi_1' \nu$$

а фазовая координата $x(t)$ определяется из следующего дифференциального уравнения

$$\dot{x}(t) = F(t) x(t) + G(t) C^{-1}(t) G(t) e^{-F(t-t_0)} \Phi_1' \nu$$

где $x(t_0)$ и $\nu$ определяются как решение системы линейных алгебраических уравнений (30). Здесь матрица $D$ формируется через (46) в следующем виде

$$D = \begin{bmatrix} 0 & -e^{F(\tau-t_0)} \Phi_2' + \Phi_1' \\ \Phi_1 - \Phi_2 e^{F(\tau-t_0)} & -\Phi_2 [W_1 - e^{F(t-t_0)} W_1 e^{F'(t-t_0)}] \Phi_2' \end{bmatrix},$$

Литература

1. Aliev F.A. Comments on 'Sweep algorithm for solving optimal control problem with multi-point boundary conditions' by M. Mutallimov, R. Zulfuqarova and L. Amirova // Adv. Differ. Equ. 2016, 131 (2016)

2. Aliev F.A. Optimization of discrete periodic systems. // Известия академии наук Азербайджанское ССР серия физико-математическич и технических наук. №3, с.32-36.

3. Aliev F.A., Larin V.B. Numeral-solution of discrete algebraic riccardi equation // Известия академии наук Азербайджанское ССР серия физико-математическич и технических наук. №5, с.94-104

4. Moszynski K. A method of solving the boundary value problem for a system of linear ordinary differensial equations // Algorytmy, 1964, v.11, №3, p.25-43.

5. Mutallimov M.M., Amirova L.I., Aliev F.A., Faradjova Sh.A., Maharramov I.A. Remarks to the paper: sweep algorithm for solving optimal control problem with multi-point boundary conditions // TWMS J. Pure Appl. Math. V.9, N.2, 2018, pp. 243-246.

6. Абрамов А.А. О переносе граничных условий для систем линейных обыкновенных дифференциальных уравнений (вариант метода прогонки). // Журнал вычислительной математики и математической физики, 1961, Т.1, №3, с. 542-545.

7. Алиев Ф.А. Задача оптимального управления линейной системой с неразделенными двухточечными краевыми условиями. Дифференциальные уравнения, 1986, № 2, с. 345-347.

8. Алиев Ф.А. Задача оптимального управления линейной системой с неразделенными двухточечными краевыми условиями. // Известия академии наук Азербайджанское ССР серия физико-математическич и технических наук. 1986, № 2, с. 345-347.

9. Алиев Ф.А. Задача оптимизации с двухточечными краевыми условиями. // Известия АН СССР, сер.техн. кибернетика, 1985, № 6, с. 138-146.

10. Алиев Ф.А. Методы решения прикладных задач оптимизации динамических систем. Баку: Елм, 1989, 320 с.

11. Алиев Ф.А. Оптимизация дискретных периодических систем в сингулярном случае // ДАН. Аз. ССР №4, с.5 1983

12. Алиев Ф.А., Исмаилов Н.А. Методы решения задач оптимизации с двухточечными краевыми условиями. Препринт / АН Азерб.ССР, Институт физики, Баку, 1985, № 151, 63 с.

13. Андреев Ю.Н. Управление конечномерными линейными объектами. Москва, Наука, 1976, 424 с.

14. Бордюг Б.А., Ларин В.Б., Тимешенко А.Г. Задачи управления шагающими аппаратами. Киев, Наук думка, 1985, 264с.

15. Брайсон А., Хо Ю Ши. Прикладная теория оптимального управления. Москва: Мир, 1972, 544с.

16. Захар-Иткин М.Х. Матричное дифференциальное уравнение Риккати и полугруппа дробно-линейных преобразований. УМН, 1973, т.28, в. 3, с. 83-120.

17. Ларин В.Б. Управление шагающими аппаратами. Киев: Наук. Думка, 1980, 168с.

18. Муталлимов М.М. Алгоритм "прогонки" для решения задачи оптимизации с неразделенными трехточечными краевыми условиями. // Докл. НАН Азербайджана, 2007, т. LXIII, № 2, с. 24-29.

19. Муталлимов М.М., Алиев Ф.А. Методы решения задач оптимизации при эксплуатации нефтяных скважин. Saarbrücken (Deutscland): LAP LAMBERT, 2012, 164 с.

20. Муталлимов М.М., Зульфугарова Р.Т., Гулиев А.П. Алгоритм прогонки для решения задач оптимального управления с неразделенными краевыми условиями. // Вестник БГУ, серия физико-математических наук, №2, 2009, Баку, с.153-160.

21. Шаманский В.Е. Методы численного решения краевых задач на ЭЦВМ. Киев: Наук.думка, ч.I, 1963, 196 с, ч.II, 1966, 244 с.

## NEW SWEEP ALGORITHM FOR SOLVING A CONTINUOUS LINEAR-QUADRATIC OPTIMIZATION PROBLEM WITH UNSEPTABLE BOUNDARY CONDITIONS

**Abstract.** A new algorithm for solving the solution of the linear-quadratic optimization problem (LQP) with unseparated boundary conditions in the continuous case is given. Using the properties of symmetry of the corresponding Hamiltonian matrix, the Euler-Lagrange equations, it is shown that linear algebraic equations for determining the missing initial data of the system being solved have a symmetric principal matrix. The results are illustrated by the example of LQP optimization (stationary case) with minimal control actions.

**Keywords:** optimization, sweep method, Riccati differensial equation, linear algebraic systems equations.